\documentclass[11pt]{article}

\begin{document}

\renewcommand{\theequation}{\arabic{section}.\arabic{equation}}
\newtheorem{lemma}{Lemma}[section]
\newtheorem{theorem}[lemma]{Theorem}
\newtheorem{remark}[lemma]{Remark}
\newtheorem{prop}[lemma]{Proposition}
 \newtheorem{coro}[lemma]{Corollary}
   \newtheorem{definition}[lemma]{Definition}
   \newcommand{\eps}{\varepsilon}
   \newcommand{\Wsob}{\smash{{\stackrel{\circ}{W}}}_2^1(D)}
   \newcommand{\EX}{{\mathbb E }}
   \newcommand{\PX}{{\mathbb P }}
\baselineskip =1.2\baselineskip
\newcommand{\be}{\begin{equation}}
   \newcommand{\ee}{\end{equation}}
   \newcommand{\tr}{\triangle}
\newcommand{\e}{\epsilon}
\renewcommand{\a}{\alpha}
\renewcommand{\b}{\beta}
\newcommand{\om}{\omega}
\newcommand{\Om}{\Omega}
\newcommand{\D}{\Delta}
\newcommand{\p}{\partial}
\newcommand{\pf}{{\noindent \sc Proof. }}
\renewcommand{\phi}{\varphi}
\newcommand{\N}{{\mathbb N}}
\newcommand{\R}{{\mathbb R}}
\newcommand{\cF}{{\cal F}}
\newcommand{\cao}{{\cal A}_1}
\newcommand{\cat}{{\cal A}_2}
\parindent0.0em

\title
{Dynamics of the Thermohaline Circulation\\ under  Wind Forcing }
\author{Hongjun Gao$^1$ and Jinqiao Duan$^2$\\
\\
1. Department of Mathematics\\Nanjing Normal University\\
Nanjing 210097,  China\\ \\
2. Department of Applied Mathematics\\ Illinois Institute of Technology \\
  Chicago, IL 60616, USA  }

\date{June 29,  2001}

\maketitle

\begin{abstract}

The ocean thermohaline circulation, also called meridional overturning circulation,
is caused by water density contrasts.  This circulation has large capacity
of carrying heat around the globe
 and it  thus affects the energy budget and further affects the climate.
We consider a thermohaline circulation model in the meridional plane
under external wind forcing.
We show that, when there is no wind forcing,
the stream function and the density fluctuation
(under appropriate metrics) tend to zero exponentially  fast as time
goes to infinity.  With rapidly oscillating wind forcing,
we obtain an averaging principle for the
thermohaline circulation model.  This averaging principle provides
convergence results and comparison estimates
between the original thermohaline circulation and the averaged
thermohaline circulation, where the wind forcing is replaced
by its time average.
This establishes the validity for using the averaged thermohaline
circulation model for numerical
simulations  at long time scales.

\vskip 0.5cm
 \noindent {\bf  Mathematics Subject Classifications}:
{Primary 35K35, 60H15, 76U05;
 Secondary 86A05, 34D35}

\vskip 0.2cm
 \noindent {\bf  Key Words}: Exponential decay,
averaging principle,
 geophysical flows, wind forcing

\end{abstract}

\section{Introduction}

In addition to the wind-driven surface circulation, the ocean also
exhibits a large meridional overturning circulation called the
thermohaline circulation. The ocean is heated (thus made less
dense) where pure freshwater is evaporated (water thus made
saltier and denser), and vice versa. The global thermohaline
circulation involves water masses sinking at high latitudes and
upwelling at lower latitudes. The process is maintained by water
density contrasts in the ocean, which themselves are created by
atmospheric forcing , namely, heat and water exchange via
evaporation and condensation.
 During the thermohaline circulation, water masses carry heat (or cold)
around the globe. Thus, it is believed that the global ocean thermohaline circulation
plays an important role in the climate \cite{Siedler}.

A two-dimensional thermohaline circulation model  involves  the
Navier-Stokes equations
 for momentum
(in the meridional plane)
 together with the convection-diffusion equations  for temperature
and for salinity.  Due to the linear equation of state (relating fluid density with
temperature and salinity),  these latter
two equations may be replaced by a single convection-diffusion equation
for density  or density fluctuation. We consider the thermohaline
circulation on a fluid domain  
in the vertical  meridional $xz$-plane, as in, e.g.,  \cite{Bouruet}:
\begin{eqnarray}
\D\psi_t + J(\D \psi, \psi ) &  = & \rho_x+ \nu \D^2 \psi   + f(x, z, t)\; ,
\label{eqn1}\\
\rho_t + J(\rho, \psi ) &  = & -N^2 \psi_x +  \frac{\nu}{Pr} \D \rho \; ,
\label{eqn2}
\end{eqnarray}
where $\psi (x,z,t)$ is the stream function which defines
 the velocity field $(u, w)= (\psi_z, -\psi_x)$;
  $\rho(x,z,t)$ is the   density   fluctuation from the
mean density; $\nu>0$ is the viscosity; $N^2>0$ is the
 mean buoyancy frequency and is taken as a constant; $Pr$ is the Prandtl number.
Finally, $J(a,b)=a_xb_y-a_yb_x$ is the Jacobian operator and
$\D = \p_{xx} + \p_{zz}$
is the Laplace operator.
Note that $\frac{\p}{\p t}+J(\cdot, \psi) = \frac{\p}{\p t}+u\p_x+w\p_z$ is the
 material derivative.
The wind forcing term $f(x, z, t)$ is to be specified below.

     In \cite{Bouruet}, the author made  some numerical simulation for
(1.1), (1.2) with periodic boundary conditions in both $x$ and $z$
for $\psi$ and $\rho$ and $f\equiv 0$.  Therefore, we assume that   $\psi$ and
$\rho$ are periodic  (with period $1$) in $x$ and $z$,  and also
assume that   $\psi$ and $\rho$
 have zero mean.

     In some recent work on the   thermohaline circulation,
the wind forcing is ignored \cite{Quon, Thual, Dijk}.
In this paper,
we will consider the impact of wind forcing on the   thermohaline circulation,
while considering the evolution of the fluid density fluctuation
(rather than  the fluid density  itself).

In the first part of this paper,  we obtain the
exponential decay estimates for the stream function
$\psi$ and density fluctuation $\rho$ with $f\equiv
0$ (no  external  wind forcing).
In the second part,  we consider   the effect of wind forcing
on the stream function and density fluctuation. We
obtain an averaging principle for rapidly oscillating wind forcing,
which provides
convergence results and comparison estimates
between the original thermohaline circulation and the averaged
thermohaline circulation.
This establishes the validity for using the averaged thermohaline
circulation model for numerical
simulations  at long time scales.


\section{Exponential Decay:    Without Wind Forcing}
\setcounter{equation}{0}

In this section, we consider the long time behavior of the
stream function and the density fluctuation in the
thermohaline circulation.
We first briefly comment on the local
existence for (1.1)-(1.2) with periodic boundary conditions 
(with period $1$) in
both $x$ and $z$ for $\psi$ and $\rho$ (with zero mean), we introduce
some notations:
$$
\int = \int_D dxdz ,
$$
where $D=\{ (x,z): 0\leq x, z \leq 1\}$ is the periodic fluid domain.
$$
H =
L_{per0}^2 = \{u: \;\; u\in L^2(D), u \mbox{ is perodic both in }
x \mbox{ and } z, \int u = 0 \},\;\; \mbox{with norm }\|\cdot\| ;
$$
$$V = H_{per0}^1 = \{u : \;\; u\in H, \nabla u \in H\},\;\;H_{per0}^2 =
\{u: \;\;  u\in H, \nabla u \in H, \tr u\in H\}, \; etc.
$$
In fact, by a result  in \cite{lady}, we know $\|\tr u\|$ is equivalent to
$\|u\|_{H_{per0}}^2$.

 Define the vorticity $\omega = \tr \psi$. It is well known that $\tr^{-1}$ exists for $\tr$ with
 periodic boundary and zero mean, then (1.1) can be written as

\be
\omega_t + J(\omega, \tr^{-1}\omega) = \rho_x + \nu \tr\omega
+ f(x, z, t). \label{2.1}
\ee

Since the nonlinear  Jacobian term  is continuous from $V\times V
\to H\times H$, by the theory of \cite{he}, we have the following
local existence result for (2.1) and (1.2):

\begin{lemma} ({\bf Local Existence})
Let $(\omega_0, \rho_0)\in V\times V$(initial values, that is
$\psi_0\in H_{per0}^3$) and $f\in L^{\infty}(0, T; H)$, then (2.1) and (1.2) with
periodic boundary conditions in both $x$ and $z$ for $\psi$ and $\rho$ with zero mean
has a unique local solution satisfying
$$\omega\in L^{\infty}(0, T; V)\cap L^{2}(0, T; H_{per0}^2)\;\;
\rho\in L^{\infty}(0, T; V)\cap L^{\infty}(0, T; H_{per0}^2),$$
that is, (1.1) and (1.2)   has a unique local solution satisfying
$$\psi\in L^{\infty}(0, T; H_{per0}^3)\cap L^{2}(0, T; H_{per0}^4)\;\;
\rho\in L^{\infty}(0, T; V)\cap L^{\infty}(0, T; H_{per0}^2),$$
where $T$ depends on $(\omega_0, \rho_0)$.
\end{lemma}

We need the following properties and estimates  (see \cite{df})
of
   the Jacobian operator $J: H_0^1\times H_0^1\to L^1$ in the sequel:
   $$
   \int_DJ(f, g)h dx dy = - \int_DJ(f, h)g dx dy, \int_DJ(f, g)g dx dy = 0,
   $$
   $$
   |\int_D J(f, g)dx dy| \le \|\nabla f\|\|\nabla g\|,
   $$
   for all $f, g, h \in H_0^1$.

 Now, we derive some a priori estimates for the solution of (1.1) and (1.2) with
 $f\equiv 0$ to ensure the global existence.  For  $f \neq 0$, the estimates are
 almost the same.

\bigskip

Multiplying (1.1) by $\psi$, performing an integration by parts
and using the periodic boundary conditions, we conclude

\be
 \frac12\frac{d}{dt}\|\nabla \psi\|^2 + \int \partial_x\rho\psi
+ \nu\|\tr\psi\|^2 = 0. \label{2.2}
\ee

Similarly, we get

\be
\frac12\frac{d}{dt}\|\rho\|^2 - N^2\int
\partial_x\rho\psi + \frac{\nu}{Pr}\|\nabla\rho\|^2 = 0.
\label{2.3}
 \ee

 Multiplying (\ref{2.2})
by $N^2$ and adding  to (\ref{2.3}),  we have
 \be
\frac12\frac{d}{dt}(N^2\|\nabla\psi\|^2 + \|\rho\|^2) +
\nu(N^2\|\tr\psi\|^2 +\frac{1}{Pr}\|\nabla\rho\|^2) = 0.
\label{2.4}
\ee

By the Poincar\'{e} inequality and taking $\alpha =
\frac{\nu}{\lambda_1}\min\{1, \frac{1}{Pr}\}$($\lambda_1$ is the
smallest eigenvalue of $-\tr$ with periodic boundary and zero
mean), we obtain
$$
 N^2\|\nabla\psi\|^2 + \|\rho\|^2 +
\int_0^t[\nu(N^2\|\tr\psi\|^2 +\frac{1}{Pr}\|\nabla\rho\|^2)] dt
$$

\be
\le e^{-\alpha t}(N^2\|\nabla\psi_0\|^2 + \|\rho_0\|^2).
\label{2.5}
\ee

If we only want to   know whether the
solution tends to zero as $t \to \infty$, we could use the
following   special Gronwall Lemma. We omit the details for this
asymptotics here, and we will concentrate on the exponential decay of
the solution in the sequel.
\begin{lemma}
 If a non-negative differential function $f$ satisfies
$$f^{\prime}(t) + \alpha_1 f(t) \le g(t),$$
where $\alpha_1 > 0$ and $\lim\limits_{t\to \infty}g(t) = 0$, then
$\lim\limits_{t\to \infty} f(t) = 0$.
\end{lemma}

\pf We first have

\be
f(t) \le e^{-\alpha_1 t}f(0) + e^{-\alpha_1 t}\int_0^t
e^{\alpha_1\tau}g(\tau) d\tau.
\ee

Since $\lim\limits_{t\to \infty}g(t) = 0$, for any given $\e > 0$,
there exists a $T$ such that
$$ g(t) \le \e, \;\; \mbox{for}\; t \ge T.$$
So, (2.6) can be written as

\be
f(t) \le e^{-\alpha_1 t}f(0) + \frac{\e}{\alpha_1}(1 - e^{-\alpha_1 t}) +
e^{-\alpha_1 t}\int_0^T e^{\alpha_1\tau}g(\tau) d\tau.
\ee

Let $t \to \infty$ and note that $\e$ is arbitrary, the result is thus obtained.

\bigskip

Now, we turn to the exponential decay estimates for the solutions.

Multiplying (\ref{eqn1}) and (\ref{eqn2}) by $- \tr\rho$ and
$\tr\psi$ respectively, integrating by parts, using the periodic
boundary conditions, we have

\be
\frac12\frac{d}{dt}\|\tr\psi\|^2 = \int\partial_x\rho\tr\psi -
\nu\|\nabla\tr\psi\|^2,
 \ee

\be
\frac12\frac{d}{dt}\|\nabla\rho\|^2 - \int J(\rho,
\psi)\tr\rho = - N^2 \int\partial_x\psi\tr\rho -
\frac{\nu}{Pr}\|\tr\rho\|^2.
\ee

Note that
\begin{eqnarray*}
\frac12\frac{d}{dt}\|\tr\psi\|^2 &=& - \int\rho\tr\psi_x -
\nu\|\nabla\tr\psi\|^2
\\
&\le& \frac{1}{4\delta\nu}\|\rho\|^2 - (1 -
\delta){\nu}\|\nabla\tr\psi\|^2
(\; \forall\; 0 < \delta < 1)\\
&\le& \frac{1}{4\delta\nu}\|\rho\|^2 - \frac{(1 -
\delta)\nu}{\lambda_1}
\|\tr\psi\|^2\\
&\le& \frac{1}{4\delta\nu}e^{-\alpha t}(N^2\|\nabla\psi_0\|^2 +
\|\rho_0\|^2)
 - \frac{(1 - \delta)\nu}{\lambda_1}\|\tr\psi\|^2.
\end{eqnarray*}
where the Poincar\'{e} inequality is used.  By  the Gronwall inequality and
(2.5), we thus obtain

$$ \|\tr\psi\|^2 \le e^{-\frac{(1 - \delta)\nu}{\lambda_1}t}\|\tr\psi_0\|^2 +
$$

\be
e^{-\frac{(1 -
\delta)\nu}{\lambda_1}t}\frac{\lambda_1}{4\delta(1-\delta)\nu^2}
\int_0^te^{-\alpha \tau}(N^2\|\nabla\psi_0\|^2 +
\|\rho_0\|^2)e^{\frac{(1 - \delta) \nu}{\lambda_1}\tau}d\tau.
 \ee

Noticing that $\alpha = \frac{\nu}{\lambda_1}\min\{1,
\frac{1}{Pr}\}$, we know for every given $Pr$, there exists a $0 <
\delta < 1$, such that

$$ \delta_1 = \alpha - \frac{(1 - \delta)\nu}{\lambda_1} > 0.$$
Hence (2.10) can be reduced to

\be
 \|\tr\psi\|^2 \le e^{-\frac{(1 -
\delta)\nu}{\lambda_1}t}(\|\tr\psi_0\|^2 +
\frac{\lambda_1}{4\delta(1-\delta)\delta_1\nu^2}(N^2\|\nabla\psi_0\|^2
+ \|\rho_0\|^2)).
\ee

Multiplying (2.8) by $ N^2$, adding to (2.9), using integrating by
parts and the fact that
 $\int[\int\partial_x\nabla\rho\nabla\psi + \int\partial_x\nabla\psi\nabla\rho] = 0$,
  we conclude that

\be
\frac12\frac{d}{dt}(N^2\|\tr\psi\|^2 + \|\nabla\rho\|^2) -
\int J(\rho, \psi)\tr\rho =  - N^2\nu\|\nabla\tr\psi\|^2 -
\frac{\nu}{Pr}\|\tr\rho\|^2.
\ee

Now we only need to estimate
$-\int J(\rho, \psi)\tr\rho = \int(\rho_z\psi_x -
\rho_x\psi_z)(\rho_{xx} + \rho_{zz})$.
Note that
\begin{eqnarray*}
\int \rho_z\psi_x\rho_{xx} &=& \frac12\int\rho_x^2\psi_{xz} - \int\rho_x\rho_z
\psi_{xx},\\
\int\rho_z\psi_x\rho_{zz} &=& - \frac12\int\rho_z^2\psi_{xz},\\
- \int\rho_x\psi_z\rho_{xx} &=& \frac12\int\rho_x^2\psi_{xz},\\
- \int\rho_x\psi_z\rho_{zz} &=& - \frac12\int\rho_z^2\psi_{xz} +
\int\rho_x\rho_z\psi_{zz}.
\end{eqnarray*}
Now we need the following lemma  about the equivalence between
 $\|u\|_{H^2}$  and
$\|\tr u\|$ for $u\in H_{per0}^2$.

\begin{lemma} For every $u\in H_{per0}^2$, we have
$$\|u\|_{H^2} \le a_1\|\tr u\|,
$$
where $a_1 = \sqrt{1 + \lambda_1 + \lambda_1^2}$.
\end{lemma}
\pf Since

\be
\|\tr u\|^2  = \int((\frac{\partial^2u}{\partial x^2})^2 + (\frac{\partial^2u}
{\partial z^2})^2 + (\frac{\partial^2u}{\partial x\partial z})^2).
\ee

Adding $\| u\|_{H^1}^2$ to both side of (2.13), and using Poincar\'{e} inequality,
 the proof of this lemma is complete.

\bigskip

Using this lemma, we imply that
\begin{eqnarray*}
|\int J(\rho, \psi)\tr\rho| &=& |\int(\rho_x^2 - \rho_z^2)\psi_{xz} - \int\rho_x
\rho_z(\psi_{xx}- \psi_{zz})|\\
&\le& \int(|\rho_x|^2 + |\rho_z|^2)|\psi_{xz}| + \int|\rho_x||\rho_z|(|\psi_{xx}| +
|\psi_{zz}|)\\
&\le& a_1(2+\sqrt{2})\|\tr\psi\|[(\int|\rho_x|^4)^{\frac12} +
\int|\rho_z|^4)^{\frac12}).
\end{eqnarray*}
By the following inequality from \cite{fmtt}
$$
\|u\|_{L^4} \le a_2\|u\|^{\frac12}\|\nabla u\|^{\frac12}, u\in
H_{per0}^1(D), 
a_2 =  (\frac1{4\pi^2}+\frac{\sqrt{2}}{\pi}+2)^{\frac14} .
$$ Therefore,
$$
 |\int J(\rho, \psi)\tr\rho| \le
a_1^2a_2^2(2+\sqrt{2})\|\tr\psi\|\|\nabla\rho\|\|\tr\rho\| $$

\be
\le \frac{\delta_2\nu}{Pr}\|\tr\rho\|^2 +
\frac{a_1^4a_2^4(6+2\sqrt{2}) Pr}{\delta_2\nu}
\|\tr\psi\|^2\|\nabla\rho\|^2,
\ee

 Combining (2.12) with (2.14) and using the Poincar\'{e} inequality, we obtain
$$ \frac12\frac{d}{dt}(N^2\|\tr\psi\|^2 + \|\nabla\rho\|^2)
$$
\be \le - N^2\nu\|\nabla\tr\psi\|^2 + (
\frac{a_1^4a_2^4(6+2\sqrt{2}) Pr}{4\delta_2\nu}\|\tr\psi\|^2 -
(1-\delta_2)\frac{\nu}{\lambda_1 Pr})\|\nabla\rho\|^2 .
\ee

 By (2.11), there exists a $t_1$ such that

$$\frac{a_1^4a_2^4(6+2\sqrt{2}) Pr}{4\delta_2\nu}\|\tr\psi\|^2 -
(1-\delta_2)\frac{\nu} {\lambda_1 Pr} = - \delta_3 < 0 ,\;0 <
\delta_3 < 1.$$

Let $\Phi_0 = \|\tr\psi_0\|^2 +
\frac{\lambda_1}{4\delta(1-\delta)\delta_1\nu^2}
(N^2\|\nabla\psi_0\|^2 + \|\rho_0\|^2)$. It follows that
$$ e^{-\frac{(1 - \delta)\nu}{\lambda_1}t}\Phi_0 <
\frac{4\delta_2(1-\delta_2)\nu^2}{a_1^4a_2^4(6+2\sqrt{2})
\lambda_1 Pr^2}.
$$

 Since $4\delta_2(1-\delta_2) \le 1$, so $t_1$ can  be chosen as

$$t_1  > \frac{\lambda_1}{(1-\delta)\nu}\ln[\frac{aa_1^4a_2^4(6+2\sqrt{2})\lambda_1
 Pr^2\Phi_0 }
{\nu^2}].$$

By (2.15), we know

 \be
 \|\nabla\rho\|^2 \le M_0,\;\; 0 \le t \le t_1
\ee

and

\be
\|\nabla\rho\|^2 \le e^{-\alpha_2 t}(N^2\|\tr\psi_0\|^2 +
\|\nabla\rho_0\|^2),\; t > t_1.
\ee

 Here
$\alpha_2 = \min\{\delta_3, \frac{N^2\nu}{\lambda_1}\}$ and $M_0$
is constant depending on $t_1$, $\|\tr\psi_0\|^2,
\|\nabla\rho\|^2, \nu, Pr$ and $\lambda_1$.

\bigskip

The estimates (2.16) and
(2.17) tell us that the mean-square norm 
of the density (fluctuation) gradient, $\|\nabla \rho\|$, 
is uniformly bounded
up to some time instant and then decay exponentially fast. 

\begin{remark} If $\nu$ is large enough or $Pr$  is small enough,
we also can have

$$\frac{a_1^4a_2^4(6+2\sqrt{2}) Pr}{4\delta_2\nu}\|\tr\psi\|^2 - (1-\delta_2)
\frac{\nu} {\lambda_1 Pr} < 0,$$

that is

$$Pr < \frac{\nu}{a_1^2a_2^2(2+\sqrt{2})}\sqrt{\frac{1}{(\lambda_1e^{-\frac{(1 - \delta)\gamma}
{\lambda_1}t}(\|\tr\psi_0\|^2 +
\frac{\lambda_1}{4\delta(1-\delta)\delta_1\nu^2}
(N^2\|\nabla\psi_0\|^2 + \|\rho_0\|^2)))}}$$

or $\nu$ satisfying
$$
e^{-\frac{(1 - \delta)\nu}{\lambda_1}t}(\|\tr\psi_0\|^2 +
\frac{\lambda_1}
{4\delta(1-\delta)\delta_1\nu^2}(N^2\|\nabla\psi_0\|^2 +
\|\rho_0\|^2))) \le 1
$$

i.e.,

$$
\nu \ge \sqrt{\frac{a_1^4a_2^4(6+2\sqrt{2})
Pr^2\lambda_1}{4\delta_2(1-\delta_2)}} .
$$
\end{remark}

\vspace{0.5cm}
   Now, we  derive the estimates for $\|\nabla\tr\psi\|$ when   the initial value
$\psi_0\in H_{per0}^3$.
   By a similar process as above, we have

\be
 \frac12\frac{d}{dt}\|\nabla\tr\psi\|^2 - \int J(\tr\psi,
\psi)\tr^2\psi =  - \int\partial_x\rho\tr^2\psi -
\nu\|\tr^2\psi\|^2.
\ee

By the same discussion  as in deriving (2.14), we have

 \be
  |\int J(v, u)\tr v| \le
a_1^2a_2^2(2+\sqrt{2})\|\tr u\|\|\nabla v\|\|\tr v\|,\;\;
\mbox{for}\; u, v\in H_{per0}^2.
 \ee

 We use inequality (2.19) for $u= \psi$ and $v= \tr\psi$, we have

$$ |\int J(\tr\psi, \psi)\tr^2\psi|  \le a_1^2a_2^2(2+\sqrt{2})\|\tr
\psi\|\|\nabla\tr\psi\|\|\tr^2\psi\| $$

\be
\le \frac{\nu}{4}\|\tr^2\psi\|^2 +
\frac{a_1^4a_2^4(6+2\sqrt{2})}{\nu}\|\tr
\psi\|^2\|\nabla\tr\psi\|^2£¬
\ee

and

\be
\int|\partial_x\rho\tr^2\psi| \le \frac{\nu}{4}\|\tr^2\psi\|^2
+ \frac{1}{\nu}\|\nabla\rho\|^2.
\ee

By (2.18), (2.20) and (2.21), we get

$$ \frac12\frac{d}{dt}\|\nabla\tr\psi\|^2 \le
\frac{1}{\nu}\|\nabla\rho\|^2 +
\frac{a_1^4a_2^4(6+2\sqrt{2})}{\nu}\|\tr
\psi\|^2\|\nabla\tr\psi\|^2 - \frac{\nu}{2}\|\tr^2\psi\|^2.
$$

Thus we conclude by using the Poincar\'{e} inequality,

 \be
\frac{d}{dt}\|\nabla\tr\psi\|^2 \le \frac{2}{\nu}\|\nabla\rho\|^2
+ (\frac{a_1^4a_2^4(12+4\sqrt{2})}{\nu}\|\tr \psi\|^2
-\frac{\lambda_1\nu}{2})\|\nabla\tr\psi\|^2.
 \ee

 Using (2.11) and for $\lambda_1$   large enough, there exists
$\alpha_3 > 0$ such
  that

$$\frac{a_1^4a_2^4(12+4\sqrt{2})}{\nu}\|\tr \psi\|^2 -
{\lambda_1\nu} = - \alpha_3 < 0,\;(\alpha_2 > \alpha_3).
$$

Hence
$$
\frac{d}{dt}\|\nabla\tr\psi\|^2 + \alpha_3\|\nabla\tr\psi\|^2 \le
\frac{2}{\nu}\|\nabla\rho\|^2,
$$
by the Gronwall's inequality, we have
$$
\|\nabla\tr\psi\|^2 \le e^{-\alpha_3t}\|\nabla\tr\psi_0 \|^2 +
\frac{2}{\nu}\int_0^t e^{-\alpha_3(t-\tau)}\|\nabla\rho\|^2 d\tau.
$$
For $t \le t_1$, using (2.16), we have

\be
\|\nabla\tr\psi\|^2 \le e^{-\alpha_3t}\|\nabla\tr\psi_0 \|^2 +
\frac{2}{\nu\alpha_3}M_0 \doteq M_1.
 \ee

 For $t > t_1$, using (2.17) and $\alpha_2 > \alpha_3$, we get

 \begin{eqnarray}
 \|\nabla\tr\psi\|^2 &\le& e^{-\alpha_3t}\|\nabla\tr\psi_0 \|^2 +
\frac{2}{\nu}(N^2\|\tr\psi_0\|^2 + \|\nabla\rho_0\|^2)\int_0^t
e^{-\alpha_3(t-\tau) - \alpha_2\tau} d\tau \nonumber\\
&=& e^{-\alpha_3t}(\|\nabla\tr\psi_0 \|^2 +
\frac{2}{\nu}(N^2\|\tr\psi_0\|^2 + \|\nabla\rho_0\|^2))\int_0^t
e^{- (\alpha_2 - \alpha_3)\tau} d\tau. \nonumber
\\&\le& e^{-\alpha_3t}(\|\nabla\tr\psi_0 \|^2 +
\frac{2}{\nu}(N^2\|\tr\psi_0\|^2 + \|\nabla\rho_0\|^2)).
\end{eqnarray}

 If $\lambda_1$ is not so large, by (2.11), there exist  $\alpha_4 > 0$ and $t_2 > 0$ large enough such that

$$\frac{a_1^4a_2^4(12+4\sqrt{2})}{\nu}\|\tr \psi\|^2 -
{\lambda_1\nu} = - \alpha_4 < 0,\;\;\mbox{ for}\; t \ge
t_2\;(\alpha_2
> \alpha_4) .
$$
Thus (2.22) can be written as

$$
\frac{d}{dt}\|\nabla\tr\psi\|^2 + \alpha_4\|\nabla\tr\psi\|^2 \le
\frac{2}{\nu}\|\nabla\rho\|^2.
$$
By the the same argument as we obtain the estimates of (2.23) and
(2.24), we could get similar estimates.

The estimates (2.23) and
(2.24) tell us that the mean-square norm 
of the vorticity gradient, $\|\nabla \omega\|$, is uniformly bounded
up to some time instant and then decay exponentially fast.

\bigskip

Thus, by the above estimates (2.5), (2.11), (2.16), (2.17), (2.23)
and (2.24), which hold  in the case of no wind forcing ($f = 0$),
we obtain the main theorem in this section:

\begin{theorem}({\bf Exponential Decay   in the Case of No Wind Forcing})
Let the initial conditions for the vorticity and
density fluctuation $(\omega_0, \rho_0)$ be in $ V\times V$
(i.e., $\psi_0\in H_{per0}^3$).
Then, when there is no external wind forcing, 
$\|\psi\|_{H_{per0}^3}$  and
  $\|\rho\|_{H_{per0}^1}$
tend to zero exponentially fast as  $t\to
\infty$.
That is,  under  appropriate norms or metrics, 
the stream function and density fluctuation
 tend to zero exponentially fast as   time goes to infinity.

Moreover,    there exists some time  instant $T>0$ such that,
the mean-square norms  for the density (fluctuation) gradient,
 $\|\nabla \rho\|$,  and  for  the vorticity gradient, $\|\nabla \omega\|$,
are uniformly bounded when $t \le T$ and exponentially decay when 
$t > T$. The  time  instant $T= \max\{t_1, t_2\}$ 
  depends on $Pr, \;N^2, \;\nu$, $\lambda_1$ and initial
values. 
\end{theorem}


\section{Averaging Principle:  Rapidly Oscillating Forcing}
\setcounter{equation}{0}

 In this section, we consider the averaging principle for the system
 of (2.1) and (1.2) under the rapidly oscillating forcing $f(x,z,t)$.
We rewrite (2.1) and (1.2) as
\begin{eqnarray} \omega_t + {\cal A}_1 \omega  &=& - J(\omega, \tr^{-1}\omega) +
\rho_x +  f(x, z, t), \label{3.1} \\
\rho_t + {\cal A}_2\rho   & = & - J(\rho, \tr^{-1}\omega ) -N^2
\tr^{-1}\omega_x  \; , \label{3.2}
\end{eqnarray}
where ${\cal A}_1$ and ${\cal A}_2$ denote the operator $- \nu
\tr$ and $- \frac{\nu}{Pr} \D$ with the periodic boundary conditions
and zero mean.
 For the rest of this section, we concentrate on
the system  (3.1)-(3.2).

We assume that  the forcing term $f$ in (3.1)  is rapidly oscillating,
i.e., it has the form
    $f(x,y, t) = f(x, y, \eta t) \doteq
   f(\eta t)$,  with parameter $\eta \gg 1.$
We also assume that $f$ has a  well-defined time average.
With such a forcing, it
is desirable to understand the fluid dynamics in some averaged
sense,  and compare the averaged flows with the original
(non-averaged) flows.

The main result of this section is an averaging principle for
(3.1)--(3.2) with rapidly oscillating forcing on finite but large
time intervals. This includes comparison estimate  and
convergence
  result  (as  $\eta \to \infty$)  between  (3.1)--(3.2) and its averaged motions.

   Starting from the fundamental
   work of Bogolyubov \cite{bo} the averaging theory for ODE has been
   developed and generalized in a large number of works (see
   \cite{bm}--\cite{fi} and the references therein).
   Bogolyubov's main theorems have been generalized in \cite{dak}
   to the case of differential equations with bounded operator-valued
   coefficients. Some problems of averaging of differential equations with
   unbounded operator-valued coefficients have been considered in
   \cite{he}--\cite{sim} in the framework of abstract parabolic equations.
   In \cite{iiy}, Ilyin considered the averaging principle for
an equation of the form

   \be
   \partial_t u = N(u) + f(\eta t),
   \ee

   where $f$ is a given forcing function and $\eta \gg 1$ is a large
   dimensionless parameter, and $f$ has a time average
 defined as

   \be
   \lim\limits_{t\to \infty}\frac{1}{t}\int_0^tf(s) ds = f_0.
   \ee

  Note that ${\cal A}_1$ and ${\cal A}_2$ are sectorial operators.
  For a sectorial operator, one can define the  fractional power of ${\cal A}$ as
follows \cite{he}:
   $${\cal A}^{\alpha} = ({\cal A}^{-\alpha})^{-1}, \;\mbox{where}\;{\cal
   A}^{-\alpha} =
   \frac{1}{\Gamma(\alpha)}\int_0^{\infty}t^{\alpha-1}e^{-{\cal A} t}  dt.$$
   The corresponding domains $D({\cal A}^{\alpha})$ are Banach spaces with
   norm given by
   $$\|x\|_{\alpha}:= \|x\|_{D({\cal A}^{\alpha})} = \|{\cal
   A}^{\alpha}x\|.$$

We recall some definitions and
  results to be used in the rest of this section.

   \begin{lemma}\cite{he} The following estimates are valid:

   \be
   \|e^{-{\cal A}t}\|_{L^2\to L^2} \le Ke^{-at}, \;\;\;\;\;\;\;\;\; t \ge
   0,
   \ee

   \be
   \|{\cal A}^{\alpha}e^{-{\cal A}t}\|_{L^2\to L^2} \le
   \frac{K_{\alpha}}{t^{\alpha}}e^{-at}, \;\; t > 0,
   \ee

   where $K, K_{\alpha}$ are positive constants.
   \end{lemma}
\begin{remark}
Since ${\cal A}_1$ and ${\cal A}_2$ are different operators but
both satisfy the conditions of {\bf Lemma 3.1}.  For the
simplicity, we take the same constants for ${\cal A}_1$ and ${\cal
A}_2$ when we use {\bf Lemma 3.1}.
\end{remark}

   \begin{lemma}\cite{he}  Given two sectorial operators $A$ and $B$ in
   $L^2$, let $D(A) = D(B)$, $Re\sigma(A) > 0, Re\sigma(B) > 0$, and for
   some $\alpha\in[0, 1)$. Let the operator $(A-B)A^{-\alpha}$ be bounded
   in $L^2$. Then for every $\gamma\in[0, 1)$, $D(A^{\gamma}) =
   D(B^{\gamma})$, the two norms being equivalent.
   \end{lemma}

 Setting
   $$ \tau = \eta t, \;\;\e = \eta^{-1},$$
   we rewrite the equations (3.1)-(3.2) in the so-called standard form

   \begin{eqnarray}
 \omega_{\tau} + \e{\cal A}_1\omega + \e J(\tr^{-1}\omega, \omega) &=& \e\rho_x +\e
f(x, y, \tau), \label{3.7} \\
\rho_{\tau} + \e{\cal A}_2\rho   + \e J(\rho, \tr^{-1}\omega )& =
&  -\e N^2\tr^{-1}\omega_x \; , \label{3.8}
\end{eqnarray}

   We assume that $f$ has a time average , $f_0(x,z)$, in $D({\cal A}^{\gamma})$;
the value of
   $\gamma$ will be specified later on.   More precisely, let $f(\tau), f_0
   \in {\cal A}^{\gamma}$ and suppose that

   \be
   \|{\cal A}^{\gamma}(\frac{1}{T}\int_t^{t+T}f(\tau) d\tau - f_0)\| \le
   \min(M_{\gamma}, \sigma_{\gamma}(T)),
   \ee

   where $M_{\gamma} > 0, \sigma_{\gamma}(T) \to 0, \;\mbox{as}\; T \to
   \infty$.

We consider the averaged equation
\begin{eqnarray}
 \bar\omega_{\tau} + \e{\cal A}_1\bar\omega + \e J(\tr^{-1}\bar\omega, \bar\omega)
  &=& \e\bar\rho_x + \e f_0(x, y), \label{3.10} \\
\bar\rho_{\tau} + \e{\cal A}_2\bar\rho   + \e J(\bar\rho,
\tr^{-1}\bar\omega )& = & -\e N^2\tr^{-1}\bar\omega_x \; ,
\label{3.11}
\end{eqnarray}
 By the method of
\cite{bav}--\cite{wang2}, we know the semigroup $S_t$
corresponding to equation (\ref{3.10})--(\ref{3.11})
    possesses absorbing sets in the space ${\bf H} = L_{per0}^2 \times L_{per0}^2,
    {\bf V} = D({\cal
   A}_1^{\frac{1}{2}}) \times D({\cal
   A}_2^{\frac{1}{2}}) = H_{per0}^1 \times H_{per0}^1$ and $D({\cal A}) = D({\cal A}_1)
   \times D({\cal A}_2)$(in fact, $D({\cal A}_1) = D({\cal
   A}_2)$).Using Lemma 3.3, we know $D({\cal A}_1^{\gamma}) = D({\cal
   A}_2^{\gamma})$, for $\gamma\in [0, 1]$.
    $\|\cdot\|$ and $\|\cdot\|_{\frac12}$ denote the norm in $L_{per0}^2$ and $H_{per0}^1$. These sets are certain
   balls $B(R_0)$ in these spaces, where $R_0$ is large enough. This means
   that for every bounded set $B$
   $$S_tB \subset B(R_0), \;\mbox{for}\; t > t_0(B, R_0).$$
 In addition, the semigroup is uniformly bounded in these spaces, that
   is, given any ball, in particular, the ball $B(R_0)$, there exists a
   ball $B(R)$ such that
   $$S_tB(R_0) \subset B(R), \;\mbox{for}\; t > 0.$$
   By increasing $R$ we may assume that
   $$S_tB(R_0) \subset B(R- r), \;\mbox{for}\; t > 0, r > 0,$$
where $r$ is a positive constant. We consider the averaging
principle in the space ${\bf V}$. Given a
   point $\omega_0$ in $B_{\bf V}(R_0)$,  we compare the  trajectories (solutions)
$(\omega(\tau), \rho(\tau)) $ and $(\bar\omega(\tau),
\bar\rho(\tau))$ of system (3.1--(3.2) and (3.7)--(3.8)
   starting from same initial point.  Consider their difference on the interval
   $\tau\in [0, \frac{T}{\e}]$, $T$ being arbitrary but fixed. We suppose
   for the moment that $(\omega(\tau), \rho(\tau)) \in B_{\bf V}(R)$. Then the difference $z(\tau) =
   \omega(\tau) - \bar\omega(\tau), \Theta = \rho - \bar\rho$ satisfies the
   equations
   \begin{eqnarray}
   \partial_{\tau} z + \e {\cal A}_1z(\tau) + \e[J(\tr^{-1}\omega, \omega) -
   J(\tr^{-1}\bar\omega, \bar\omega)] &=& \e((\rho_x - \bar\rho_x) + (f(\tau) - f_0)), \\
   \partial_{\tau}\Theta + \e {\cal A}_2(\tau) + \e[J(\tr^{-1}\omega, \rho) -
   J(\tr^{-1}\bar\omega, \bar\rho)] &=& -\e N^2(\tr^{-1}\omega_x - \tr^{-1}
   \bar\omega_x).
\label{eqnz}
   \end{eqnarray}

   We first give the following lemma, the proof can be obtained by
   direct estimate.
\begin{lemma}
   The nonlinear operator $J(u, v)$ is a bounded Lipschitz map in the
   following sense:
   $$
   \|J(u_1, v_1) - J(u_2, v_2)\| \le
   $$
   \be
   C_{\frac{1}{2}}( \|u_1\|_{\frac{1}{2}} + \|u_2\|_{\frac{1}{2}} +
   \|v_1\|_{\frac{1}{2}} + \|v_2\|_{\frac{1}{2}})(\|u_1 -
   v_1\|_{\frac{1}{2}} + \|u_2 - v_2\|_{\frac{1}{2}}),
   \ee

where $C_{\frac{1}{2}}$ is some positive constants.
   \end{lemma}

   Inverting the linear operators ${\cal A}_1$ and ${\cal A}_2$ we
   come to the equivalent integral equations of (3.12) and (3.13)
\begin{eqnarray}
   z(\tau) &=& - \e\int_0^{\tau}e^{-\e{\cal A}_1(\tau -s)}[J(\tr^{-1}\omega,
   \omega) - J(\tr^{-1}\bar\omega, \bar\omega)]ds \nonumber\\&+&
    \e\int_0^{\tau}e^{-\e{\cal A}_1(\tau -s)}(\rho_x - \bar\rho_x)ds +
    \e\int_0^\tau e^{-\e{\cal A}_1(\tau -s)}(f(s) - f_0)ds,
   \\
\Theta (\tau)&=& - \e\int_0^{\tau}e^{-\e{\cal A}_2(\tau
-s)}[J(\tr^{-1}\omega, \rho) -
   J(\tr^{-1}\bar\omega, \bar\rho)]ds\nonumber\\  &-&
   \e N^2\int_0^{\tau}e^{-\e{\cal A}_2(\tau
-s)}(\tr^{-1}\omega_x - \tr^{-1}
   \bar\omega_x).
  \end{eqnarray}
  Using (3.6) and (3.14),  we see that the $\|\cdot\|_{\frac12}$-norm of the first
   term in the right hand side of (3.15) satisfies the inequality
      $$\|\e\int_0^{\tau}\cao^{\frac{1}{2}}e^{-\e{\cao}{(\tau
   -s)}}[J(\tr^{-1}\omega, \omega) - J(\tr^{-1}\bar\omega,
   \bar\omega)]ds\|
   $$
   $$\le \e\int_0^{\tau}K_{\frac12}C_{\frac12}\e^{-\frac12}(\tau -s)^{-\frac12}e^{-\e
   a(\tau - s)}2R\|z(s)\|_{\frac12}ds
   $$

\be
   = 2RK_{\frac12}C_{\frac12}\e^{\frac12}\int_0^{\tau}(\tau -s)^{-\frac12}e^{-\e
   a(\tau - s)}\|z(s)\|_{\frac12}ds.
   \ee

   Let us estimate the second term in the right hand side of (3.15).

\be
  \|\e\int_0^{\tau}\cao e^{-\e{\cal A}_1(\tau -s)}(\rho_x -
   \bar\rho_x)ds\| \le K_{\frac12}\e^{\frac12}\int_0^{\tau}(\tau -s)^{-\frac12}
   e^{-\e a(\tau - s)}\|\Theta(s)\|_{\frac12}ds.
   \ee

  Now let us estimate the third term in the right hand side of (3.15).
Integrating by parts we have
$$
   \|\e\int_0^\tau e^{-\e{\cal A}_1(\tau -s)}(f(s) - f_0)ds\|_{\frac12}
   $$
   $$
   = \|-\e e^{-\e{\cal A}_1(\tau -s)}\int_s^{\tau}(f(t) - f_0)dt|_0^{\tau} +
   \e^2\int_0^\tau \cao e^{-\e{\cal A}_1(\tau -s)}\int_s^{\tau}(f(s) -
   f_0)ds\|_{\frac12}
   $$
   $$
   \le
   \|\e\cao^{\frac12-\gamma}e^{-\e{\cal
   A}_1\tau}{\cal A}_1^{\gamma}\int_0^{\tau}(f(t) - f_0)dt\|
   $$

   \be
    + \|\e^2\int_0^\tau \cao^{\frac32-\gamma} e^{-\e{\cal A}_1(\tau
   -s)}\cao^{\gamma}\int_s^{\tau}(f(s) - f_0)ds\|.
   \ee

   Using (3.6) and (3.9), we  further have
   $$
   \|\e\cao^{\frac12-\gamma}e^{-\e{\cal
   A}_1\tau}{\cal A}_1^{\gamma}\int_0^{\tau}(f(t) - f_0)dt\| \le
   \e K_{\frac12-\gamma}e^{-\e a\tau}(\e\tau)^{\gamma -
    \frac12}\|\frac{1}{\tau}\int_0^{\tau}{\cal A}_1^{\gamma}(f(t) - f_0)dt\|
   $$
   $$
   = (\e\tau)^{\frac12 + \gamma}K_{\frac12-\gamma}e^{-\e
   a\tau}\|\frac{1}{\tau}\int_0^{\tau}{\cal A}_1^{\gamma}(f(t) - f_0)dt\|
   $$

   \be
   \le (\e\tau)^{\frac12 + \gamma}K_{\frac12-\gamma}\min(M_{\gamma},
   \sigma_{\gamma}(\tau))e^{-\e a\tau} =: L(\tau).
   \ee

   For any $\delta > 0$, let $\tau_{\delta}$ be so large that for $\tau \ge
   \tau_{\delta},  we have \sigma_{\gamma} \le \delta$.
Let $\e_0$ be so small that
   for $\e < \e_0$ the  inequality $\frac{T}{\e} > \tau_{\delta}$ is valid.
   Then
   $$
   L(\tau) \le G_{\gamma 1}(T, \e) = e^{-\e a\tau}\left\{\begin{array}{cc}
   T^{\frac12 + \gamma}K_{\frac12 - \gamma}\delta, \;&\mbox{if} \;\tau \ge
   \tau_{\delta},\\
   (\e\tau)^{\frac12 + \gamma}K_{\frac12-\gamma}M_{\gamma}, \;&\mbox{if}\;
   \tau < \tau_{\delta}.
   \end{array}
   \right.
   $$
   Let $\gamma > - \frac12$.  Note that $\tau_{\delta}$ does not depend on $\e$.
   Taking $\delta\to 0$ and then $\e\to 0$, we obtain

   \be
   \|\e e^{-\e{\cal A}_1\tau}\int_0^{\tau}(f(t) - f_0)dt\|_{\frac12} \le
   G_{\gamma 1}(T, \e) \to 0\;\mbox{when}\; \e \to 0.
   \ee

   $$
   \|\e^2\int_0^\tau \cao^{\frac32-\gamma} e^{-\e{\cal A}_1(\tau
   -s)}\cao^{\gamma}\int_s^{\tau}(f(s) - f_0)ds\| \le K_{\frac32 -
   \gamma}\e^{\frac12 + \gamma}\int_0^{\tau}\min(M_{\gamma},
   \sigma_{\gamma}(u))u^{\gamma - \frac12}du
   $$
   $$
   \le K_{\frac32 - \gamma}M_{\gamma}\e^{\frac12 +
   \gamma}\int_0^{\tau_{\mu}}u^{\gamma - \frac12}du + K_{\frac32 -
   \gamma}\e^{\frac12 + \gamma}\mu\int_0^{\frac{T}{\e}}u^{\gamma - \frac12}du
   $$

   \be
   = K_{\frac32 - \gamma}(\gamma + \frac12)^{-1}(M_{\gamma}(\e\tau_{\mu})^{\frac12 +
   \gamma} + \mu T^{\frac12 + \gamma}) =: G_{\gamma 2}(T, \e),
   \ee

   where for any $\mu > 0$ we have chosen $\tau_{\mu}$ so large that
   $\sigma_{\gamma}(\tau) < \mu$ when $\tau > \tau_{\mu}$. Letting $\mu \to
   0$ and then $\e\to 0$ we obtain
   $$ G_{\gamma 2}(T, \e) \to 0, \;\; \e\to 0.$$
   Thus, by (3.17)--(3.22) we obtain the following inequality:

   \be
   \|z(\tau)\|_{\frac12} \le K\e^{\frac12}\int_0^{\tau}(\tau
   -s)^{-\frac12}(\|z(s)\|_{\frac12} + \|\Theta(s)\|_{\frac12})ds +  G_{\gamma}(T, \e),
   \ee

   where $K = \max\{2RK_{\frac12}C_{\frac12}, K_{\frac12}\}$ and $G_{\gamma} = G_{\gamma 1} + G_{\gamma 2} \to
   0, \e \to 0$.

   Similar to the argument  above for $z(\tau)$, and using (3.6) and
   (3.14), we get

  \be
   \|\Theta\|_{\frac12} \le K_1\e^{\frac12}\int_0^{\tau}(\tau
   -s)^{-\frac12}(\|z(s)\|_{\frac12} + \|\Theta(s)\|_{\frac12})ds,
   \ee

   where $ K_1$ depends on $R, K_{\frac12}, C_{\frac12}$ and $\lambda_1$.

Adding (3.23) and (3.24), we finally obtain

  \be
  \|z(\tau)\|_{\frac12} + \|\Theta(\tau)\|_{\frac12} \le (K + K_1)\e^{\frac12}\int_0^{\tau}(\tau
   -s)^{-\frac12}(\|z(s)\|_{\frac12} + \|\Theta(s)\|_{\frac12})ds + G_{\gamma}(T,
   \e).
   \ee

\bigskip

  Here we need the following fact.
   \begin{lemma}\cite{he} Let $\gamma\in(0, 1]$ and for $t\in [0, T]$
   $$u(t) \le a + b\int_0^t(t-s)^{\gamma-1}u(s)ds.$$
   Then
   $$u(t) \le aE_{\gamma}((b\Gamma(\gamma))^{\frac{1}{\gamma}}t),$$
   where the function $E_{\gamma}(z)$ is monotone increasing and
   $E_{\gamma}(z)\sim \gamma^{-1}e^z$ as $z\to \infty$.
   \end{lemma}

\bigskip

    Applying this lemma to the inequality (3.25) on $\tau\in[0,
   \frac{T}{\e}]$, we obtain
   \begin{eqnarray}
   \|z(t)\|_{\frac12} + \|\Theta(\tau)\|_{\frac12} &\le&\nonumber\\
   G_{\gamma}(T, \e)
   E_{\frac12}(\e\tau\pi (K + K_1)^2) &\le&
   G_{\gamma}(T, \e)E_{\frac12}(T\pi (K + K_1)^2)  :=   \eta_T^1(\e).
   \end{eqnarray}
   We thus have proved the proximity of solutions of (3.1) and (3.2)
 in ${\bf V}$,  for
   the  trajectory  $(\omega(t), \rho(t))$   with  initial condition
$(\omega(0), \rho(0)) \in B_{\bf V}(R_0)$ stays in the ball $B(R)$
on the  time interval $[0, \frac{T}{\e}]$.

    Let $\e$ be so small that the right-hand side of (3.26) are
   less than $\frac{r}{2}$, where $r$ is defined earlier in this section
when we discuss absorbing sets. Suppose that the trajectory
$(\omega(t), \rho(t))$
   leaves the ball $B(R)$ during the interval $[0, \frac{T}{\e}]$ and let
   $\tau^*$ be the first moment where
$\|\omega(\tau^*)\|_{\frac12} + \|\rho(\tau^*)\|_{\frac12} = R$.
However, on
    the interval $\tau\in[0, \tau^*]$ both trajectories stay in the ball
   $B(R)$ and what we have proved so far shows that the inequality
   $$
\|\omega(\tau) - \bar\omega(\tau)\|_{\frac12} 
+ \|\rho(\tau) - \bar\rho(\tau)\|_{\frac12}
   \le \frac{r}{2}
$$
is valid. In
   particular, it is valid for $\tau = \tau^*$. This together with the
   inequality
$\|\bar\omega(\tau^*)\|_{\frac12} + \|\rho(\tau^*)\|_{\frac12}  \le R - r$,
   which holds by the
   hypothesis of the following theorem and the property of the semigroup
   $S(t)$, gives the contradiction
   $$\|\omega(\tau^*)\|_{\frac12} +  \|\rho(\tau^*)\|_{\frac12} \le
   \|\omega(\tau^*) - \bar\omega(\tau^*)\|_{\frac12} +  \|\rho(\tau^*) -
   \bar\rho(\tau^*)\|_{\frac12}
     + \|\bar\omega(\tau^*)\|_{\frac12}
     +  \| \bar\rho(\tau^*)\|_{\frac12} \le R - \frac{r}{2},
   $$
   since $\|\omega(\tau^*)\|_{\frac12} + \|\rho(\tau^*)\|_{\frac12} =
   R$.

Thus we have the main theorem in this section:

\begin{theorem}
 ({\bf Averaging Principle in the Case of
Rapidly Oscillating Wind Forcing}) Assume that  the wind forcing
has a time average. Let $T > 0$ be arbitrary and fixed. 
If $\gamma > - \frac12$ and    the initial conditions 
for the vorticity and density fluctuation 
$(\omega(0), \rho(0)) = (\bar\omega(0), \bar \rho(0))$ are in
the absorbing ball  $B_{\bf V}(R_0)$ (depending on $\gamma$),  
then for $\tau\in[0, \frac{T}{\e}]$, we have the following 
comparison and convergence estimate
between the thermohaline circulation and the averaged thermohaline
 circulation 
$$
 \|\omega(\tau) - \bar\omega(\tau) \|_{\frac12} +
 \| \rho(\tau) - \bar\rho(\tau) \|_{\frac12} \le
    \eta_T^1(\e)\to 0, \mbox{as}\;\e\to 0,
$$
where $\eta_T^1(\e)$ is a decaying function defined in (3.26).
\end{theorem}

\section{Summary}

The ocean thermohaline circulation has important impacts on the climate.
We have considered a thermohaline circulation model in the meridional plane
under external wind forcing.
We have shown that, when there is no wind forcing,
the stream function and the density fluctuation
(under appropriate metrics) tend to zero exponentially  fast as time
goes to infinity (Theorem 2.5).  With rapidly oscillating wind forcing,
we have obtained an averaging principle for the
thermohaline circulation model (Theorem 3.6).  This averaging principle provides
convergence results and comparison estimates
between the original thermohaline circulation and the averaged
thermohaline circulation, where the wind forcing is replaced
by its time average.


\bigskip

{\bf Acknowledgements.}   A part of this work was done at the
Oberwolfach Mathematical Research Institute, Germany and Institute
of Mathematics and Its Applications,  while J. Duan was a Research
in Pairs Fellow, supported by {\em Volkswagen Stiftung}. This work
was partly supported by the NSF Grant DMS-9973204 and and by the
grant of NNSF of China 10001018.   And a part of this work was done
while H. Gao was visiting   Illinois Institute of Technology, Chicago, and
Institute of Mathematics and Its Applications, Minnesota, USA.



\begin{thebibliography}{99}

\bibitem{Bouruet} P. Bouruet-Aubertot, C. Koudella, C. Staquet and K. B. Winters,
Particle dispersion and mixing induced by breaking internal gravity waves,
{\em Dynamics of Atmos. Oceans} {\bf 33} (2001), 95-134.

\bibitem{df} V. P. Dymnikov and A. N. Filatov,
{\em Mathematics of
 Climate Modeling},
Birkhauser, Boston, Cambridge, MA, 1997.

\bibitem{Constantin} P. Constantin and C. Foias,
{\em Navier-Stokes Equations}, Univ. of Chicago Press, Chicago,
1988.

\bibitem{lady} O. A. Ladyzhenskaya, {\em The Boundary Value Problems of
Mathematics Physics}, Springer-Verlag, 1985.

\bibitem{fmtt} C. Foias, O. Manley, R. Temam and Y. M. Treve,
Asymptotic analysis of the Navier-Stokes equations,
Phys. D, 9(1983), 157-188.

\bibitem{bo} N. N. Bogolyubov, On some statistical methods in
   mathematical physics, Izdat. Akad. Nauk Ukr. SSR, Kiev 1945.
   \bibitem{bm} N. N. Bogolyubov and Yu. A. Mitropolskii, {\em Asymptotic
   methods in the theory of non-linear oscillations}, English transl.,
   Gordon and Breach, New York, 1962.

   \bibitem{mi} Yu. A. Mitropolskii, {\em The methods of averaging in
   non-linear mechanics}, Naukova Dumka, Kiev 1971(Russian).

   \bibitem{fi} A. N. Filatov, {\em Asymptotic methods in the theory of
   differential and integrodifferential equations}, Fan, Tashkent,
   1974 (Russian).

   \bibitem{dak} Y. L. Daletskii and M. G. Krein, {\em Stability of
   solutions of differential equations in Banach space}, English transl.,
   Amer. Math. Soc., Providence, RI 1974.

   \bibitem{he} D. Henry, {\em Geometric theory of semilinear parabolic
   equations}, Springer-Verlag, New York, 1981.

   \bibitem{lez} B. M.  Levitan and V. V. Zhilov, Almost periodic functions
   and differential equations, English transl., Cambridge Univ. Press,
   Cambridge, 1982.

\bibitem{Verhulst} F. Verhulst,
On averaging methods for partial differential equations,
{\em preprint}, 1999.

   \bibitem{sim} I. B. Simonenko, Justification of the method of averaging
   for abstract parabolic equations, English transl. in Math. USSR-Sb.
   10(1970)53--61.

   \bibitem{iiy} A. A. Ilyin, Averaging principle for dissipative dynamical
   system with rapidly oscillating right-hand sides, Math. Sb., 187(1996),
   635--677.

\bibitem{bav} A. V. Babin and M. I. Vishik, {\em Attractor of evolution
   equations}, English transl., North-Holland, Amsterdam, 1992.

   \bibitem{hale} J. K. Hale, {\em Asymptotic behavior for dissipative
   dynamical system},
   Amer. Math. Soc., Providence, RI, 1988.

\bibitem{Marotzke} J. Marotzke, Abrupt climate change and thermohalie
circulation: Mechanisms and predictability, {\em Proc. National Acad. Sci.},
{\bf 97} (2000), 1347-1350.



 \bibitem{wang2} S. Wang, Attractors for the 3D baroclinic
 quasi-geostrophic equations of large scale atmosphere, J. Math.
 Anal. Appl., 165(1992), 266-283.

\bibitem{Dijk}H. A. Dijkstra and J. D. Neelin,
Imperfections of the thermohaline circulation: Latitudinal asymmetry and preferred
northern sinking, {\em J. Climate} {\bf 13} (2000), 366-382.

\bibitem{DuanSchm} J. Duan and B. Schmalfu{\ss},
The 3D Quasigeostrophic Fluid Dynamics under Random
Forcing on Boundary,   submitted,  2000.

\bibitem{Gill} A. E. Gill, {\em Atmosphere-Ocean Dynamics},
Academic Press, New York, 1982. (Chapter 9)

\bibitem{Ler98}M.~Leroux, {\em Dynamic Analysis of Weather and Climate}.
 John Wiley \& Sons, 1998. (Chapter 11)


\bibitem{Oz} T. Ozgokmen and E. P. Chassignet,
A numerical study of two-dimensional
turbulent gravity currents descending a slope, preprint, 2001.

\bibitem{Ped87}
J.~Pedlosky,
 {\em Geophysical Fluid Dynamics},
Springer Verlag, New-York, Berlin, 1987.

\bibitem{PeiOor92}J.~P. Peixoto and A.~H. Oort, {\em Physics of Climate},
 Springer, New York, 1992.

 \bibitem{Quon} C. Quon and M. Ghil,
Multiple equilibria in thermosolutal convection
 due to salt-flux boundary conditions,
{\em J. Fluid Mech.} {\bf 245} (1992), 449-483.

\bibitem{Siedler} G. Siedler, J. Church and J. Gould,
{\em Ocean Circulation and Climate: Observing and Modelling the
Global Ocean}, Academic Press, San Diego, USA, 2001.

\bibitem{Thual} O. Thual and J. C. McWilliams, The catastrophe structure of
thermohaline convection in a two-dimensional fluid model and a comparison with
low-order box model, {\em Geophys. Astrophys. Fluid Dynamics} {\bf 64} (1992),
 67-95.

\bibitem{Wright} D. G. Wright  and T. F. Stocker,
A zonally averaged  ocean  model for
 the thermohaline circulation. Part I: Model development and flow dynamics,
{\em J. Phys. Oceanography} {\bf 21} (1991),  1713-1724.

\bibitem{Washington} W. M. Washington and C. L. Parkinson,
{\em An Introduction to
Three-Dimensional Climate Modeling}, Oxford Univ. Press, 1986.

\end{thebibliography}
\end{document}